\documentclass{article}[fleqn]

\usepackage{arxiv}

\usepackage[utf8]{inputenc} 
\usepackage[T1]{fontenc}    
\usepackage{hyperref}       
\usepackage{url}            
\usepackage{booktabs}       
\usepackage{amsfonts}       
\usepackage{nicefrac}       
\usepackage{microtype}      
\usepackage{lipsum}
\usepackage{graphicx}

\usepackage{array}
\usepackage{multirow}
\usepackage{amssymb,amsmath,amsfonts,etoolbox}
\usepackage[usenames, dvipsnames]{color}
\usepackage[footnotesize]{subfigure}
\usepackage{caption}
\usepackage{subfigure}
\usepackage{microtype}
\usepackage{epstopdf}
\usepackage{hyperref}
\usepackage{color}
\usepackage{subfigure}
\usepackage[inline]{enumitem}
\usepackage{booktabs}
\usepackage{amsfonts}
\usepackage{amssymb}
\usepackage{amsmath, amsthm}
\usepackage{bm}
\usepackage[ruled,vlined,linesnumbered]{algorithm2e}

\usepackage[normalem]{ulem}
\usepackage{mathtools}

\usepackage{pgf}
\usepackage{arydshln}
\usepackage{tikz}
\usetikzlibrary{arrows,automata,positioning}
\usepackage{datatool}
\usepackage{booktabs}
\usepackage{natbib}
\usepackage{csvsimple}
\bibpunct[, ]{(}{)}{,}{a}{}{,}%
\def\newblock{\ }%
\usepackage{graphicx,adjustbox}
\usepackage[normalem]{ulem}
\usepackage[linesnumbered]{algorithm2e}
\usepackage{algpseudocode}
\usepackage{bbm}
\usepackage{tcolorbox}
\usepackage{float}
\usepackage{wrapfig}
\usepackage{lscape}
\usepackage{rotating}
\usepackage{epstopdf}

\usepackage{amsthm}
\DeclareMathOperator*{\argminB}{argmin}   

\title{A ride time-oriented scheduling algorithm for dial-a-ride problems}

\author{
 Claudia Bongiovanni \\
 Polytechnique Montr\'{e}al and CIRRELT\\
\texttt{claudia.bongiovanni@polymtl.ca}
\And
Nikolas Geroliminis \\
\'{E}cole Polytechnique F\'{e}d\'{e}rale de Lausanne\\
\texttt{nikolas.geroliminis@epfl.ch}
\And
Mor Kaspi \\
Tel-Aviv University \\
\texttt{morkaspi@tauex.tau.ac.il}
}

\begin{document}
\maketitle
\begin{abstract}
This paper offers a new algorithm to efficiently optimize scheduling decisions for dial-a-ride problems (DARPs), including problem variants considering electric and autonomous vehicles (e-ADARPs). The scheduling heuristic, based on linear programming theory, aims at finding minimal user ride time schedules in polynomial time. The algorithm can either return optimal feasible routes or it can return incorrect infeasibility declarations, on which feasibility can be recovered through a specifically-designed heuristic. The algorithm is furthermore supplemented by a battery management algorithm that can be used to determine charging decisions for electric and autonomous vehicle fleets. Timing solutions from the proposed scheduling algorithm are obtained on millions of routes extracted from DARP and e-ADARP benchmark instances. They are compared to those obtained from a linear program, as well as to popular scheduling procedures from the DARP literature. Results show that the proposed procedure outperforms state-of-the-art scheduling algorithms, both in terms of compute-efficiency and solution quality.
\end{abstract}

\keywords{scheduling \and routing \and battery management \and dial-a-ride \and heuristics}

\section{Introduction}
    Many real-world urban mobility and supply chain problems involve time-dependent allocation and sequencing decisions that must be frequently re-optimized over time. The re-optimization process is typically composed of the production and evaluation of a plethora of potential solutions. The evaluation of these solutions is performed by scheduling algorithms, whose goals are to determine the timing of the sequential decisions, while respecting several time-related constraints and objectives \citep[e.g.,][]{vidal2015timing}. Their efficient resolution is particularly relevant for large-scale and dynamic problems, since in these cases scheduling decisions must be made frequently and quickly, as for example in the context of large neighborhood search methods \citep[e.g.,][]{gschwind2019adaptive}. In addition to computational efficiency, it is essential that scheduling algorithms provide high-quality solutions, given their direct impact on the quality of allocation decisions and thus on the overall performance of the optimization strategy.
    
    Scheduling algorithms are of utmost importance for hardly-constrained vehicle routing problems. They assume feasible vehicle routes in terms of capacity and precedence constraints and are designed to provide feasible decisions on time-related variables (e.g., service start times and waiting times at each node). Demonstrating the feasibility of a fixed route and composing high-quality vehicle schedules is a challenging task for most vehicle routing problems because of the presence of time-related constraints (e.g., time windows and maximum route duration). For dial-a-ride (DARP) problems, this task is made even more challenging by the presence of constraints on maximum user ride time \citep[][]{cordeaulaporte2007}, as well as constraints on battery management and charging time for fleets composed of electric and autonomous vehicles \citep[e-ADARP,][]{bongiovanni2019,su2022column, su2023deterministic}. 
    Consequently, the goal of some scheduling algorithms for DARPs is optimize time dependent level of service while satisfying several time related constraints. The former is often measured by the user excess ride time, that is the delay users experience for sharing rides as compared to a taxi service \citep[e.g., see Section 3 in][]{molenbruch2017}.

    A variety of heuristics have been proposed in the DARP literature to manage path planning. However, these procedures typically face a trade-off between returning correct feasibility declarations and increasing the quality of the returned solutions. \cite{cordeaulaporte2003tabu} propose an 8-step procedure setting the earliest start time at each vertex in the route and using forward slack times to delay the start times at pickup locations in vision of maximum ride time constraints \citep{savelsbergh1992}. \cite{parragh2009heuristic} observe that adopting a sequential approach to avoid maximum ride time violations does not necessarily minimize the total user excess ride time in the schedules. In fact, delaying the service start time at a given pickup location may decrease the excess ride time of the specific request but increase the excess ride time for other requests in the route. As such, \cite{parragh2009heuristic} modify the procedure in \cite{cordeaulaporte2003tabu} by adapting the computation of forward slack times such that increases in the user excess ride time of any request in the route is avoided. As a result, the returned feasible schedules minimize the total user excess ride time at the expense of incorrect infeasibility declarations. This last aspect is a drawback that is tackled in \cite{molenbruch2017multidir}, who propose a procedure starting by considering a possibly travel-time infeasible schedule setting the excess ride time of each user at its lowest bound. Infeasibility relates to travel time shortages between successive nodes and is succesively recovered by shifting service start times such that the total user excess ride time is minimized.

	In this work, we propose a novel scheduling heuristic which is numerically shown to provide excess-time optimal solutions in most cases, although it may result with incorrect infeasibility declarations or suboptimal solutions. The designed algorithm is complemented by a battery management heuristic to deal with DARP extensions employing electric autonomous vehicles, i.e., e-ADARPs. The goal of the battery management heuristic is to assign charging times at the visited stations. Indeed, a schedule minimizing excess ride time may not be battery feasible, and vice versa. The procedure builds on two main observations: \begin{enumerate*}[label=(\roman*)] \item  minimizing excess-time while respecting time-window and user-ride-time constraints is in fact equivalent to assigning the right amount of waiting time at all nodes in the routes; and \item ensuring battery feasibility is in fact possible by recharging as much as possible at the visited facilities, as early as possible.\end{enumerate*} The proposed scheduling and battery management algorithms are tested on millions of DARP and e-ADARP routing solutions obtained by developing an adaptive large neighborhood search heuristic \citep[][]{ropke2006adaptive}. The quality of the obtained scheduling solutions are compared against: \begin{enumerate*}[label=(\roman*)] \item a linear program; \item the well-known 8-step scheduling procedures by \cite{cordeaulaporte2003tabu} and \cite{parragh2009heuristic}; and \item the procedure by \cite{molenbruch2017multidir}. \end{enumerate*} 

	The rest of the paper is organized as follows: Section~\ref{problem_def} introduces the excess-ride-time scheduling problem, Section~\ref{route_eval} provides the scheduling algorithm, Section~\ref{battery} presents the battery management problem and the heuristic procedure to provide battery-feasibility. Finally, Section~\ref{experiments} provides numerical experiments comparing the scheduling algorithm against state-of-the-art procedures, and Section~\ref{conclusion} summarizes the main concepts of this paper and provides an overlook to future research.

	\section{Scheduling Problem}\label{problem_def}

	Consider a predetermined route sequence $\mathcal{I}$ of $M$ nodes, which includes pickup locations, dropoff locations, and potentially some charging stations. Without loss of generality, assume that the sequence satisfies routing constraints, as well as precedence and load constraints. Note that the given sequence may not necessarily satisfy all timing and battery management constraints. Then, the optimization problem consists in scheduling the service start times in the sequence as to minimize the total user excess ride time, guarantee battery, time window, and maximum ride time feasibility. 

	Start by considering a specific sub-sequence $\mathcal{\bar{I}}$ of $\bar{M}\leq M$ nodes, which is one of the derived sub-sequences obtained by splitting $\mathcal{I}$ by the visited charging stations. Without loss of generality, assume that the visited charging station at the beginning of $\mathcal{\bar{I}}$ represents an origin depot and the charging station at the end of $\mathcal{\bar{I}}$ represents a destination depot. Furthermore, denote by $i \in \{1,\ldots,n\}$ the set of requests contained in sequence $\mathcal{\bar{I}}$, $P_i$ the set of pickups, and $D_i$ the set of dropoffs. Pickup locations are characterized by loads $l_{P_i}$ and dropoff locations by loads $l_{D_i}=-l_{P_i}$. Furthermore, all locations $i,j \in \mathcal{\bar{I}}$ feature service times $s_i$, direct travel times $t_{i,j}$, and time windows [$arr_i$, $dep_i$] limiting the time at which service may start. Note that, it is possible to set time windows around the pickup locations and derive the time windows around the corresponding dropoff locations by means of the maximum ride times $u_{P_i}$, and vice versa \citep[][]{cordeaulaporte2003tabu}. The goal is to be able to optimize $\mathcal{\bar{I}}$ with respect to the total excess ride time by setting optimal service start times $T_i$ with $i\in\{1,\ldots,\bar{M}\}$ in vision of maximum ride time and time window constraints. As such, the scheduling problem for sub-sequence $\mathcal{\bar{I}}$ can be stated as the following linear program:\\
 
	\small
	\setlength{\belowdisplayskip}{0pt} \setlength{\belowdisplayshortskip}{0pt}
	\setlength{\abovedisplayskip}{0pt} \setlength{\abovedisplayshortskip}{0pt}
	\begin{equation} \label{eq:1b}
	(LP1) \hspace{0.2cm} \min \sum\limits_{i\in \{1,\ldots,n\}} (T_{D_i}-T_{P_i}-d_{P_i}-t_{P_i,D_i})
	\end{equation}
	\begin{equation*}
	s.t.:
	\end{equation*}
	\begin{equation} \label{eq:2b}
	T_i + t_{i,i+1} + s_i \leq T_{i+1}  \hspace*{0.5cm}   \forall i \in \{1,2,\ldots,\bar{M}-1\}
	\end{equation}
	\begin{equation} \label{eq:3b}
	T_{D_i} - T_{P_i} - d_{P_i} \leq u_{P_i}  \hspace*{0.5cm}   \forall i \in \{1,\ldots,n\}
	\end{equation}
	\begin{equation} \label{eq:3c}
	arr_i \leq T_i \leq dep_i \hspace*{0.5cm}   \forall i \in  \{1,2,\ldots,\bar{M}\}
	\end{equation}\\
	\normalsize
	where the objective function~(\ref{eq:1b}) minimizes the the sum of excess ride times of all requests in $\mathcal{\bar{I}}$, constraints~(\ref{eq:2b}) set the service start times between consecutive nodes, while constraints (\ref{eq:3b})-(\ref{eq:3c})
	impose maximum ride time and time window constraints respectively.

	Note that time windows [$arr_i$, $dep_i$] can be tightened in light of the travel times and service times between consecutive nodes. As such, by sequentially inspecting the sequence, it is possible to calculate the earliest time $ET_i$ and latest time $LT_i$ at which service can start at node $i$ by using the following recursive formulas:\\
	\small
	\setlength{\belowdisplayskip}{0pt} \setlength{\belowdisplayshortskip}{0pt}
	\setlength{\abovedisplayskip}{0pt} \setlength{\abovedisplayshortskip}{0pt}
	\begin{equation} \label{eq:3d}
	ET_i=\max\{arr_i, ET_{i-1}+ t_{i-1,i}+s_i\} \hspace*{0.5cm} \forall i \in \{2,\ldots,\bar{M}\}, ET_1=arr_1
	\end{equation}
	\begin{equation} \label{eq:3e}
	LT_i=\max \{ dep_i, LT_{i+1} - t_{i,i+1}-s_i \} \hspace*{0.5cm} \forall i \in \{1,\ldots,\bar{M}-1\}, LT_{\bar{M}}=dep_{\bar{M}}
	\end{equation}\\
    \normalsize
	Hence, a tighter formulation to LP1 can be obtained by substituting constraints~(\ref{eq:3c}) with:\\
	\small
	\setlength{\belowdisplayskip}{0pt} \setlength{\belowdisplayshortskip}{0pt}
	\setlength{\abovedisplayskip}{0pt} \setlength{\abovedisplayshortskip}{0pt}
	\begin{equation} \label{eq:3f}
	ET_i\leq T_i\leq LT_i \hspace*{0.5cm} \forall i \in  \{1,2,\ldots,\bar{M}\}
	\end{equation}\\
	\normalsize
	Service start time $T_i$ at node $i$ depends on the initial departure time from the depot, the total travel time up to $i$, the total service time spent serving nodes before $i$, and the total vehicle waiting time up to $i$. That is:\\
	\small
	\setlength{\belowdisplayskip}{0pt} \setlength{\belowdisplayshortskip}{0pt}
	\setlength{\abovedisplayskip}{0pt} \setlength{\abovedisplayshortskip}{0pt}
	\begin{equation}\label{eq:1d}
	T_i=T_1+ \sum\limits_{j=1}^{i-1} t_{j,j+1} + \sum\limits_{j=1}^{i-1} d_{j} + \sum\limits_{j=1}^{i} w_j
	\end{equation}\\
	\normalsize
	Here, $w_i$ denotes the waiting time at node $i$. Note that, service may start as soon as possible without penalizing the objective function, i.e. $T_1=ET_1$. With such representation of service start times, the objective function (\ref{eq:1b}) can be re-written as follows:\\
 	\small
	\setlength{\belowdisplayskip}{0pt} \setlength{\belowdisplayshortskip}{0pt}
	\setlength{\abovedisplayskip}{0pt} \setlength{\abovedisplayshortskip}{0pt}
	\begin{equation*} \label{eq:3h}
	\min \sum\limits_{i\in \{1,\ldots,n\}} (T_1+ \sum\limits_{j=1}^{D_i-1} t_{j,j+1} + \sum\limits_{j=1}^{D_i-1} d_{j} + \sum\limits_{j=1}^{D_i} w_j) - (T_1+ \sum\limits_{j=1}^{P_i-1} t_{j,j+1} + \sum\limits_{j=1}^{P_i-1} d_{j} + \sum\limits_{j=1}^{P_i} w_j) - d_{P_i} -t_{P_i,D_i} =
	\end{equation*}
	\begin{equation} \label{eq:3i}
	\min \sum\limits_{i\in \{1,\ldots,n\}} (\sum\limits_{j=P_i}^{D_i-1} t_{j,j+1} + \sum\limits_{j=P_i}^{D_i-1} d_{j} + \sum\limits_{j=P_i}^{D_i} w_j - d_{P_i} - t_{P_i,D_i})
	\end{equation}\\
	\normalsize
	Since travel times between consecutive nodes and the service times are deterministic parameters, minimizing the objective function~(\ref{eq:3i}) is equivalent to:\\
    \small
	\setlength{\belowdisplayskip}{0pt} \setlength{\belowdisplayshortskip}{0pt}
	\setlength{\abovedisplayskip}{0pt} \setlength{\abovedisplayshortskip}{0pt}
	\begin{equation} \label{eq:3l}
	\min \sum\limits_{i\in \{1,\ldots,n\}} \sum\limits_{j=P_i+1}^{D_i} w_j = \min \sum\limits_{i=1}^{\bar{M}}L_i w_i
	\end{equation}\\
	\normalsize
	Where $L_i=\sum\limits_{j=1}^{i-1}l_j$ represents the vehicle load up to node $i$.

	Equivalent to the re-writing of the objective function, constraints (\ref{eq:3f}) can be re-defined through (\ref{eq:1d}) as follows: \\
	\small
	\setlength{\belowdisplayskip}{0pt} \setlength{\belowdisplayshortskip}{0pt}
	\setlength{\abovedisplayskip}{0pt} \setlength{\abovedisplayshortskip}{0pt}
	\begin{equation}\label{eq:2dd}
	\sum\limits_{j=1}^{i} w_j \geq ET_i - \sum\limits_{j=1}^{i-1} t_{j,j+1} + \sum\limits_{j=1}^{i-1} d_j - ET_1 \hspace*{0.5cm} \forall i \in \{1,2,\ldots,\bar{M}\}
	\end{equation}
	\begin{equation}\label{eq:3dd}
	\sum\limits_{j=1}^{i} w_j \leq LT_i - \sum\limits_{j=1}^{i-1} t_{j,j+1} + \sum\limits_{j=1}^{i-1} d_j - ET_1 \hspace*{0.5cm} \forall i \in \{1,2,\ldots,\bar{M}\}
	\end{equation}\\
	\normalsize
	Note that these constraints provide lower and upper bounds to the total waiting time amount that needs to be distributed between $i$ and all nodes preceding $i$. As such, the linear program (LP1) can be equivalently re-defined as the following linear program:\\
	\small
	\setlength{\belowdisplayskip}{0pt} \setlength{\belowdisplayshortskip}{0pt}
	\setlength{\abovedisplayskip}{0pt} \setlength{\abovedisplayshortskip}{0pt}
	\begin{equation} \label{eq:1c}
	(LP2) \hspace{0.2cm} \min \sum\limits_{i=1}^{\bar{M}} L_iw_i
	\end{equation}	
	\begin{equation*}
	s.t.:
	\end{equation*}
	\begin{equation} \label{eq:2dd}
	\sum\limits_{j=1}^{i} w_j \geq ET_i - \sum\limits_{j=1}^{i-1} t_{j,j+1} - \sum\limits_{j=1}^{i-1} d_j - ET_1 \hspace*{0.5cm} \forall i \in \{1,2,\ldots,\bar{M}\}
	\end{equation}
	\begin{equation} \label{eq:3dd}
	\sum\limits_{j=1}^{i} w_j \leq LT_i - \sum\limits_{j=1}^{i-1} t_{j,j+1} - \sum\limits_{j=1}^{i-1} d_j - ET_1 \hspace*{0.5cm} \forall i \in \{1,2,\ldots,\bar{M}\}
	\end{equation}
	\begin{equation} \label{eq:3ee}
	\sum\limits_{j=i+1}^{D_i} w_j \leq u_i - \sum\limits_{j=i}^{D_i -1} t_{j,j+1} - \sum\limits_{j=i+1}^{D_i -1} d_j \hspace*{0.5cm} \forall i \in \mathcal{P}
	\end{equation}\\
	\normalsize
	Remark that constraints~(\ref{eq:2b}) from (LP1) are guaranteed by the definition of equation~(\ref{eq:1d}), which now composes constraints~(\ref{eq:2dd}) and ~(\ref{eq:3dd}). The right-hand side of (\ref{eq:2dd}) and (\ref{eq:3dd}) represent the minimal total waiting time that must be assigned up to node $i$, and the maximal total waiting time that can be assigned up to node $i$, without violating the time windows at sucessive nodes in the sequence. For convenience, denote the right-hand side of constraints (\ref{eq:2dd}) and (\ref{eq:3dd}) by $\Delta_i$ and $\Theta_i$ respectively. Using equations (\ref{eq:3d}) and (\ref{eq:3e}) we have that $\Delta_i \leq \Delta_{i+1}\hspace*{0.2cm} \forall i \in \{1,2,\ldots,M-1\}$ and $\Theta_i \leq \Theta_{i+1}\hspace*{0.2cm} \forall i \in \{1,2,\ldots,M-1\}$. That is, the total minimal and maximal waiting time that needs to be distributed in the sequence may only increase between consecutive nodes. As such, note that the total minimal waiting time $\Delta_M$, i.e. at the destination depot, represents a waiting time amount that cannot be avoided in any feasible solution. Finally, the objective of the problem reduces to optimally distribute $\Delta_{M}$ among all nodes $\{1,\ldots,M\}$ in consideration of the total load $L_i$ at each node in the sequence, which impacts the total user excess ride time.

	\section{Scheduling Procedure}\label{route_eval}
		
	In order to solve (LP2), we propose the procedure reported in the pseudo-code from Algorithm~(\ref{alg:algorithm1}) and explained next. The algorithm proceeds in the sense of the sequence and checks that, for each encountered node $i \in \{1,\ldots, \bar{M}\}$, a minimal total waiting time $\Delta_i$ has been assigned up to node $i$. If the algorithm detects a total waiting time shortage at node $i$, i.e. $\sum\limits_{k=1}^{i} w_k < \Delta_i$, the total waiting time at $i$ and its preceding nodes needs to be increased. Given that the objective function depends on the total vehicle load $L_i$, the algorithm starts by considering adding waiting time to nodes $j\leq i$ featuring the lowest minimal total load up to $i$, i.e. $j=\argminB_{k\in\{1,\ldots,i\}}L_k$. Note that, if multiple nodes featuring an equivalent minimal total load exist,	the first node can be selected without loss of generality. For node $j$, the total waiting time can be feasibly increased by a maximum amount $\delta w_j$ defined by constraints (\ref{eq:3dd}) and (\ref{eq:3ee}). That is, while deciding upon an increment in waiting time at node $j$, one needs to check that excess-ride time constraints are not violated for requests whose pickups precede $j$ and whose dropoffs follow $j$. Furthermore, in order to guarantee time-window feasibility of the whole sequence, one needs to check that an increment in waiting time at $j$ does not exceed the maximal waiting time that can be assigned to up to node $j$, i.e. $\Theta_j - \sum\limits_{k=1}^{j} w_k$. Finally, waiting time at node $j$ can be incremented by $\delta w_j$, which is computed as the minimum between the amount defined by excess-ride time constraints, time-window constraints, and the total waiting time shortage defined by $\Delta_i -\sum\limits_{k=1}^{i} w_k$. After the update of $w_j$, if node $j$ has reached its maximum waiting time limit by updating $\Delta_i$, that is $\delta w_j< \Delta_i -\sum\limits_{k=1}^{i} w_k$, node $j$ is removed from the list $\Omega$ of potential nodes whose waiting time may be further increased.
	If $\sum\limits_{k=1}^{i} w_k<\Delta_i$, i.e. there is still a total waiting time shortage at $i$, the total waiting time is increased at the next node $\bar{j}$ up to $i$ featuring the second lowest total vehicle load. This iterative process terminates as soon as $\sum\limits_{k=1}^{i} w_k \geq \Delta_i$, that is when sufficient waiting time has been assigned to $i$ and all of its preceding nodes. In this case, the algorithm moves inspecting $i+1$ and up to the end of the sequence. If at the end of the whole process, $\Delta_{\bar{M}}$ has been feasibly assigned, the algorithm terminates with a feasible solution. Section~\ref{incorrect_infeas} shows that if a feasible solution is not obtained, this may lead to an incorrect infeasibility declaration due to maximum ride time violations. On these solutions, it is possible to recover feasibility through the recourse heurisitc presented in Section~\ref{recourse_heur}, which is applied at step 17. in Algorithm~(\ref{alg:algorithm1}). Note that the procedure between line 1. and 15. of Algorithm~(\ref{alg:algorithm1}) results in a worst-time complexity of $O(\bar{M}^2)$, since in the worst case, the step reported in line 6. may be executed $\bar{M}$ times and the procedure in line 8. contains at most $\bar{M}$ components. However, the application of the recourse heuristic leads to a worst-time complexity of $O(n \times \bar{M}^2)$, since,
    Algorithm~(\ref{alg:algorithm1}) may be restarted as many times as the number of requests $n$. Finally, note that Algorithm~(\ref{alg:algorithm1}) can either terminate with a feasible solution, a correct infeasibility declaration, a suboptimal solution, or with an incorrect infeasible solution, on which feasibility can be heuristically recovered. 

 	\RestyleAlgo{algoruled}
	\begin{algorithm}[H]
		\scriptsize
		\KwIn{vehicle route sequence ($\bar{I}=\{1,\ldots,\bar{M}\}$), pickups $P_i$, dropoffs $D_i$, earliest start times $ET_i $, latest start times $LT_i$, maximum ride times $u_{P_i}$, service durations $s_i$, travel times $t_{i,j}$}
		\KwOut{Waiting times $w_i$ with $i \in \bar{I}$, feasibility $check$}
		initialize $w_i.=0\hspace{0.2cm} \forall i \in \{1,\ldots,\bar{M}\}$\;
		initialize $\Omega=\emptyset$ \;
		initialize $start$ $node=1$\;
		initialize $check=true$\;
	    \While{$check=true$}{
		\For{$i \rightarrow$ $start$ $node$:$\bar{M}$}{
			Update: $\Omega=\Omega \cup \{i\}$ \;
			\While{$\sum\limits_{k=1}^{i} w_k < \Delta_i$}{
				Set: $j=\argminB_{k\in\Omega} L_k$ \;
				Compute: $\delta w_j = \min\{\min_{k\in \{P|k\leq j \hspace*{0.1cm}\&\hspace*{0.1cm} n+k\geq j\}} u_k - \sum\limits_{l=k}^{(n+k)-1} t_{l,l+1} - \sum\limits_{l=k+1}^{(n+k)-1} d_l - \sum\limits_{l=k+1}^{n+k} w_l;\hspace{0.1cm}\Theta_j - \sum\limits_{k=1}^{j} w_k;\hspace{0.1cm}\Delta_i - \sum\limits_{k=1}^{i} w_k\}$\;
	            Set: $ w_j = w_j + \delta w_j$\;
	            \If{$\delta w_j< \Delta_i -\sum\limits_{k=1}^{i} w_k$}{Update: $\Omega=\Omega \setminus \{j\}$\;}
	            \If{$\sum\limits_{k=1}^{i} w_k\geq\Delta_i$}{\text{break}\;}
				\If{$\sum\limits_{k=1}^{i} w_k<\Delta_i$ $\&$ $\Omega=\emptyset$}{Employ: $Recourse$ $heuristic$ (Algorithm~\ref{alg:algorithm_rec})\; \If{Recourse heuristic is not successful}{$check=false$\; \text{break}\;}\Else{$start$ $node=j+1 \leq i$\; Update waiting times\; Restart algorithm from line 6.}}}}}
		\caption{Ride time-oriented scheduling algorithm}
		\label{alg:algorithm1}
	\end{algorithm}

	\subsection{Incorrect infeasibility declarations}\label{incorrect_infeas}
	    
	   There may be cases in which, at the end of Algorithm~(\ref{alg:algorithm1}), the total waiting time shortage $\Delta_{\bar{M}}$ cannot be feasibly assigned. In particular, it is possible that the waiting times at all nodes preceding $i$ have been updated, i.e., $\Omega=\emptyset$, but there is still a shortage of waiting time at node $i$, i.e., $\sum\limits_{k=1}^{i} w_k<\Delta_i$. In this case, the Algorithm~(\ref{alg:algorithm1}) prematurely terminates and the sequence is deemed infeasible. However, note that, due to the myopic nature of the procedure, this may result in an incorrect infeasibility declaration. In fact, the procedure is designed to assign waiting times to nodes $j \leq i$ based on the waiting times shortage $\Delta_i$ and without considering the waiting times shortages at subsequent nodes $\{i+1, \ldots,\bar{M}\}$. To ensure feasibility at subsequent nodes after detecting a waiting times shortage $\Delta_i$, it may be necessary to apply waiting times to nodes $j \leq i$ that: \begin{enumerate*} \item are not the first node among nodes featuring the lowest total load; and \item do not feature the lowest total load \end{enumerate*}. Failing to apply waiting times at the right nodes without knowledge of future waiting times shortages may result in situations in which the nodes in $\Omega$ cannot be further pushed forward in time without violating the maximum ride time constraints of nodes in $\Omega$ and, consequently, an incorrect infeasibility declaration. 

    	\begin{figure}[tb]
    	\centering
    	\resizebox{14cm}{!} {
    	\includegraphics{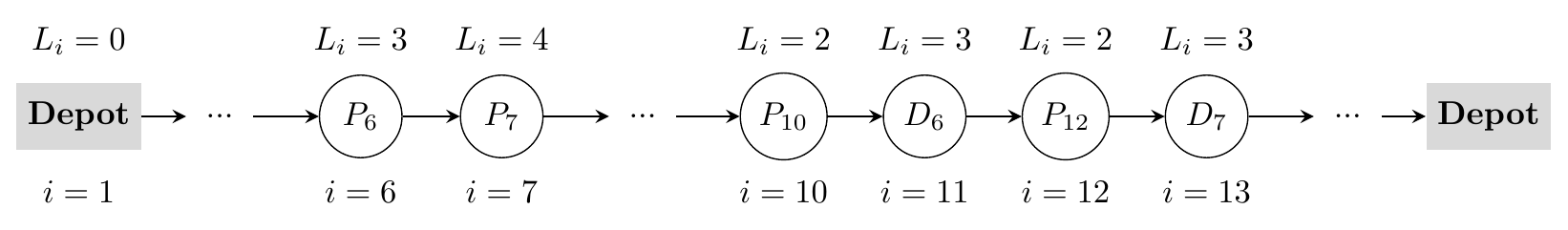}}
    	\caption{Example of a route leading to an incorrect infeasibility declaration}
    	\label{fig:example}
        \end{figure} 
	   
	   This drawback is illustrated in the routing sequence in Figure~\ref{fig:example} and explained next. Suppose that, at iteration ten, Algorithm~(\ref{alg:algorithm1}) detects a waiting time shortage $\Delta_{P_{10}}>0$ and that all of this waiting time shortage is applied at pickup node $P_{10}$, given that it features the lowest total load among all of its preceding nodes. As a consequence of increasing the waiting times at node $P_{10}$, suppose that the service start time at node $D_6$ reaches the latest service start time $dep_{D_6}$ and that the ride time of node $P_7$, with dropoff location $D_7$, reaches its maximum. The algorithm proceeds inspecting the sequence and a new waiting times shortage is detected at node $P_{12}$, i.e., $\Delta_{P_{12}}>0$.  However, it is not possible to increase waiting times at nodes preceding $P_{12}$, given the time window constraints at node $D_6$, and it is not possible to increase the waiting times at node $\Delta_{P_{12}}$, given the maximum ride time constraints at node $D_7$. In this case, the scheduling algorithm prematurely terminates and returns an infeasibility declaration for the given sequence. However, this is an incorrect infeasibility declaration given that a feasible solution can be constructed by applying the waiting time shortage $\Delta_{P_{10}}$ at node $P_{7}$ instead of node $P_{10}$, although the total load at node $P_{7}$ is higher than the total load at node $P_{10}$. This reduces the ride time of node $P_7$, with dropoff location $D_7$ following node $P_{12}$.  Consequently, the waiting times shortage $\Delta_{P_{12}}$ can be applied at node $P_{12}$, without violating the maximum ride time constraints of node $P_7$ and a feasible solution is returned. To prevent incorrect infeasibility declarations, we propose a recourse heuristic, whose aim is to revisit previously-made waiting time decisions and decrease the chances of obtaining incorrect infeasibility declarations. The details of the recourse heuristic are provided in Algorithm~\ref{alg:algorithm_rec} and are explained next. 
    
    \subsection{Recourse heuristic to reduce incorrect infeasibility declarations}\label{recourse_heur}
    
    \RestyleAlgo{algoruled}
	\begin{algorithm}[tb]
		\scriptsize
		\KwIn{vehicle route sequence ($\bar{I}=\{1,\ldots,\bar{M}\}$), pickups $P_i$, dropoffs $D_i$, earliest start times $ET_i $, latest start times $LT_i$, maximum ride times $u_{P_i}$, service durations $s_i$, travel times $t_{i,j}$, waiting times $w_i$, remaining waiting time shortening $\Delta_i$}
		\KwOut{start node $start$ $node$, waiting times $w_i$ with $i \in \bar{I}$, feasibility $check$}
		initialize $start$ $node=1$\;
		initialize $check=false$\;
		find all nodes $j_l\leq i$, $l \in \{1,\ldots,\tilde{M}\}$ that have reached maximum ride time constraints\; 
		\For{$l \rightarrow$ 1:$\tilde{M}$}{
			Set: $\tilde{w}_i=w_i$, $\tilde{w}_i=0$ $\forall i \in \{j_{l+1}, \ldots, i\}$\; 
			Recompute: $LT_i$ \;
			Compute: $\delta \tilde{w}_{j_l}$ as in step 10. of  Algorithm~\ref{alg:algorithm1}\;
			\If{$\delta w_{j_l}>0$}{
	        Update: $\tilde{w}_{j_l} = \tilde{w}_{j_l} + \delta w_j$\;
	        Update: $start$ $node=j_l$\; 
	        Update: $check=true$\; 
	        Set:  $w_i=\tilde{w}_i$\;
	        break\;
	        }
			}
		\caption{Recourse algorithm}
		\label{alg:algorithm_rec}
	\end{algorithm}
	
	Suppose that the scheduling algorithm has detected a waiting times shortage $\Delta_i$ at node $i$ and that, after inspecting all nodes $j\leq i$, $\Delta_i>0$ and $\Omega=\emptyset$. In this case, it may be necessary to modify waiting time decisions that were taken at steps preceding $i$ through a heuristic procedure, e.g.,  Algorithm~\ref{alg:algorithm_rec}. In particular, there may exists nodes that reached their maximum ride time as a consequence of waiting time decisions at previous nodes featuring the lowest total loads. Nodes $j\leq i$ that reached their maximum ride time can be detected by employing the first term in step 10. of Algorithm~\ref{alg:algorithm1}. Suppose that $\tilde{M}$ nodes $j_1<j_2<\ldots<j_i$ have reached their maximum ride time at iteration $i$. This means that, in order to recover a feasible route, the sum of the waiting time at these nodes need to be shifted by a total of $\Delta_i$. Note that postponing the service start time at $j_l$, $l \in \{1,\ldots,\tilde{M}\}$, automatically postpones the service start times at any successive nodes, including nodes $j_2, \ldots, j_i$. Furthermore, it may modify waiting time decisions at nodes following $j_l$. As such, the recourse heuristic iteratively explores possibilities to postpone service start times at nodes $j_l$, $l \in \{1,\ldots,\tilde{M}\}$, by employing the conditions defined by constraints (\ref{eq:3dd}) and (\ref{eq:3ee}) and employed in step 10. of Algorithm~(\ref{alg:algorithm1}). Note that the proposed recourse algorithm is a heuristic, given that the waiting time shortage could be applied at nodes $j_l$ as well as nodes preceding $j_l$. If the waiting times at node $j_l$ can be postponed by a maximum amount $\delta w_{j_l} >0$, the waiting time at $j_l$ is updated, the waiting times at all nodes following $j_l$ are re-initialized, and Algorithm~(\ref{alg:algorithm1}) is re-started from $j_{l+1}$. Otherwise, if $\delta w_{j_l} =0$ $\forall j_l$, $l \in \{1,\ldots,\tilde{M}\}$, then the solution is deemed infeasible and Algorithm~(\ref{alg:algorithm1}) is prematurely terminated. 
	
	\section{Battery Management Problem}\label{battery}

	In the case of electric and autonomous vehicles, i.e., an e-ADARP, it is necessary to show that the retrieved excess-time optimal schedules are also battery-feasible. Indeed, Algorithm~(\ref{alg:algorithm1}) disregards battery considerations which are instead part of e-ADARPs. Battery-feasibility aspects can be implemented by integrating battery-related decision variables and constraints into (LP1), as presented in Section 2.2 in \cite{bongiovanni2019}, and summarized next.

    Consider the predetermined route sequence $\mathcal{I}$ of $M$ nodes. As described in Section~\ref{problem_def}, this route originates at an origin depot, terminates at a destination depot, visits a total of $n$ requests, and may pass through $N$ charging stations $s \in \mathcal{S} \subset \mathcal{I}$. The latter are necessary to refill the battery level $B_i$ of the vehicle after spending battery inventory at rate $\beta_{i,i+1}$, for $i\in \{1,\ldots,M-1\}$. The vehicle can only access charging stations $s \in \mathcal{S}$ when no requests are onboard and is allowed recharge for a variable duration $E_s$ at rate $\alpha_s$. Furthermore, it departs from the origin depot with a pre-defined vehicle battery inventory $B_{init}$ and is required to reach the destination depot with a minimal vehicle battery inventory $B_{end}$. The latter is defined as a proportion of the total vehicle battery capacity $Q$, by means of a minimal battery level ratio $r$. Under these settings, the integrated scheduling and battery management problem is represented as (LP1), supplemented with the following constraints:\\
   	\small
	\setlength{\belowdisplayskip}{0pt} \setlength{\belowdisplayshortskip}{0pt}
	\setlength{\abovedisplayskip}{0pt} \setlength{\abovedisplayshortskip}{0pt} 
    \begin{equation}\label{eq:batt1}
    B_{i} = B_{init}  \hspace*{0.5cm} i=1
    \end{equation}
    \begin{equation}\label{eq:batt2}
    B_{i+1} \leq B_{i} - \beta_{i,i+1} \hspace*{0.5cm} \forall i\in \{1,\ldots,M-1\} \setminus \mathcal{S}, i \neq j
    \end{equation}
    \begin{equation}\label{eq:batt3}
    B_{i+1} \geq B_{i} - \beta_{i,i+1} \hspace*{0.5cm} \forall i\in \{1,\ldots,M-1\} \setminus \mathcal{S}
    \end{equation}
    \begin{equation}\label{eq:batt4}
    B_{s+1} \leq B_{s} + \alpha_sE_s - \beta_{s,s+1} \hspace*{0.5cm} s \in \mathcal{S}
    \end{equation}
    \begin{equation}\label{eq:batt5}
    B_{s+1} \geq B_{s} + \alpha_sE_s - \beta_{s,s+1} \hspace*{0.5cm} s \in \mathcal{S}
    \end{equation}
    \begin{equation}\label{eq:batt6}
    Q \geq B_{s} +\alpha_{s}E_{s} \hspace*{0.5cm} \forall s \in \mathcal{S}
    \end{equation}   
    \begin{equation}\label{eq:batt7}
    B_{i} \geq rQ \hspace*{0.5cm} i=M
    \end{equation}  
    \begin{equation} \label{eq:batt8}
    E_{s} \leq T_{s+1} - t_{s,s+1} - T_{s} \hspace*{0.5cm} \forall s \in \mathcal{S}
    \end{equation}
    \begin{equation} \label{eq:batt9}
    E_{s} \geq T_{s+1} - t_{s,s+1} - T_{s} \hspace*{0.5cm} \forall s \in \mathcal{S}
    \end{equation}
    \begin{equation}\label{eq:batt10}
    B_{i} \geq 0 \hspace*{0.5cm} \forall i \in \mathcal{I}
    \end{equation}
    \begin{equation}\label{eq:batt11}
    E_{s} \geq 0 \hspace*{0.5cm} \forall s \in \mathcal{S}
    \end{equation}
    
    \subsection{Battery management heuristic} 
    	\RestyleAlgo{algoruled}
	\begin{algorithm}[tb]
	   \scriptsize
	   \KwIn{vehicle route sequence ($I=\{1,\ldots,\bar{M}\}$), Waiting times $w_i$ ($i \in \mathcal{I}$), charging facilities $\{s_1,\ldots,s_N\}$, charging rates $\alpha_s$, discharging rate $\beta$, travel times $t_{i,i+1}$, earliest start times $ET_i$, latest start times $LT_i $, State of charge $B_i$}
	   \KwOut{Charging times $E_s$ with $s \in \{s_1,\ldots,s_N\}$, Boolean battery feasibility check}
	   initialize $check=true$\;
	   initialize $B_1=B_0$ \;
	   \While{$check=true$}{
	      \For{$i \rightarrow$ 2:$\bar{M}$}{
	           Compute: $B_i=B_{i-1}-\beta \times t_{i,j}$\;
			   \If{$B_{i}<0$}{the sub-sequence is infeasible $\rightarrow$ $check=false$\;}
	           \If{$i \in \mathcal{S}$}{Set: $E_i=\min\{ LT_{i+1} - ET_i; (Q-B_{i})/\alpha_{i}\}$\; Set: $B_{i}=B_{i}+E_{i}$\; Update: $w_i$ $\forall i \in \{i,\ldots,\bar{M}\}$\;}
			   }
	     }
	\caption{Battery management heuristic}
	\label{alg:algorithm2}
	\end{algorithm}
	The battery management problem, defined by (LP1) supplemented with constraints (\ref{eq:batt1})-(\ref{eq:batt11}), is solved through the battery management heuristic in Algorithm~\ref{alg:algorithm2}, which builds on the assumption that battery levels can be seen as an inventory which can only decrease with traveling. Then, for feasibility purposes, it is always better to recharge as much as possible, as early as possible. If the schedule that maximizes battery recharging does not satisfy the imposed battery management constraints, no other recharging plan may be feasible, i.e., the given sequence and excess-time optimal schedule is declared infeasible. This is due to the fact that the scheduling algorithm presented in Algorithm~\ref{alg:algorithm1} is designed to minimize the completion time and, consequently, to provide the shortest schedules and the longest charging opportunities. Consequently, if battery feasible solutions do exist, at least one of them minimizes the total user excess ride time. Note that the proposed recharging procedure in Algorithm~\ref{alg:algorithm2} is a heuristic since it cannot be formally proven optimal (e.g. maximizing the charging time at station $i$ may decrease opportunities to recharge more at following charging stations).
	
    The battery management heuristic in Algorithm~\ref{alg:algorithm2} considers the route schedule computed by employing Algorithm~\ref{alg:algorithm1} and sequentially computes the battery inventory between successive nodes in consideration of the initial battery level $B_1$, the battery discharge rate $\beta$, and the total travel time up to $i$. Note that battery discharge can be equally computed by energy consumption models \citep{goeke2015,pelletier2017}. During this iterative process, when a charging facility is encountered, its maximal recharging time is bound by:
	\begin{enumerate*}[label=(\roman*)] 
	\item the difference between the service start time at the current node and the latest service start time at the following node,
	and \item the time needed to fully recharge.
	\end{enumerate*}
	To ensure feasibility, the recharging time at the station needs to be set to the minimum of the two. After the charging time is calculated, the station's battery level and the waiting times of the station and subsequent nodes are updated. The procedure prematurely terminates only if a node $i$ features a negative battery inventory, after which the route is declared battery-infeasible.
	
	\section{Numerical Results}\label{experiments}
	
    Numerical experiments are performed on the DARP instances presented in \cite{cordeaulaporte2003tabu} and on the e-ADARP instances presented in \cite{bongiovanni2019}. Specifically, we extract routing solutions by employing the adaptive large neighborhood search described in \cite{ropke2006adaptive}, with 1,000 iterations. For each problem instance, the initial solution is constructed from the best known solution, if available, or by employing a  construction heuristic in which a random subset of customers are assigned to the vehicles, and served in a taxi mode (i.e., pickup-dropoffs). At every iteration of the large neighborhood search, we solve the scheduling problem through a linear program and compare its results to the route evaluation procedures by \cite{cordeaulaporte2003tabu}, \cite{parragh2009heuristic}, and \cite{molenbruch2017multidir}, as well as this paper. This results in a thorough comparison of the scheduling algorithms on 21,844,304 feasible DARP solutions and 2,370,243 feasible e-ADARP solutions. The numerical experiments are implemented in Julia v1.8.2 and run on 2 x AMD Rome 7532 @ 2.40 GHz 256M cache L3 CPU clusters. Each instance is run on five of such CPUs with 186 Gb of RAM. The LP is implemented in the JuMP modeling language \citep{dunning2017jump} v1.3.1 and solved with Gurobi v 0.11.3. 

     \begin{filecontents*}{excel/darp/Table_logs.csv}
scenario,number of routes,number feasible routes,range of route nodes,average cpu LP [ms],BONG: number incorrect feasible declarations,BONG: number incorrect infeasible declarations,BONG: number infeasible declarations because of TW violations,BONG: number infeasible declarations because of RT violations,BONG: number infeasible declarations because of BATT violations,BONG: number of deviating solutions,BONG: average relative deviation,BONG: average cpu algorithm [ms],BONG: total number of recourse,BONG: avg recourse time,BONG: avg route completion time,BONG: total route wait time,CORD: number incorrect feasible declarations,CORD: number of incorrect infeasible declarations,CORD: number infeasible declarations because of TW violations,CORD: number infeasible declarations because of RT violations,CORD: number infeasible declarations because of BATT violations,CORD: number of deviating solutions,CORD: average relative deviation,CORD: average cpu algorithm [ms],CORD: total number of recourse,CORD: avg recourse time,CORD: avg route completion time,CORD: total route wait time,PARR: number incorrect feasible declarations,PARR: number incorrect infeasible declarations,PARR: number infeasible declarations because of TW violations,PARR: number infeasible declarations because of RT violations,PARR: number infeasible declarations because of BATT violations,PARR: number of deviating solutions,PARR: average relative deviation,PARR: average cpu algorithm [ms],PARR: total number of recourse,PARR: avg recourse time,PARR: avg route completion time,PARR: total route wait time,MOL: number incorrect feasible declarations,MOL: number incorrect infeasible declarations,MOL: number infeasible declarations because of TW violations,MOL: number infeasible declarations because of RT violations,MOL: number infeasible declarations because of BATT violations,MOL: number of deviating solutions,MOL: average relative deviation,MOL: average cpu algorithm [ms],MOL: total number of recourse,MOL: avg recourse time,MOL: avg route completion time,MOL: total route wait time
pr01,122148,122148,4 to 24,1.960,0,0,0,0,0,0,0.000,0.649,0,0.000,362.658,234.885,0,0,0,0,0,32553,19.679,0.462,0,0,355.986,234.885,0,0,0,0,0,13,6.915,0.439,0,0,357.553,110.583,0,0,0,0,0,3,5.349,0.360,0,0,352.670,234.885
pr02,228021,228021,4 to 32,2.058,0,0,0,0,0,0,0.000,0.661,0,0.000,405.022,221.451,0,0,0,0,0,63355,18.742,0.627,0,0,396.197,221.451,0,7,0,7,0,858,4.802,0.621,0,0,400.803,136.039,0,7,0,7,0,150,3.923,0.419,0,0,395.056,221.451
pr03,447151,447151,4 to 34,2.059,0,0,0,0,0,0,0.000,0.852,0,0.000,403.094,178.997,0,0,0,0,0,121082,14.201,0.719,0,0,395.181,178.997,0,11,0,11,0,883,6.533,0.716,0,0,396.376,85.074,0,11,0,11,0,250,4.673,0.441,0,0,393.150,178.996
pr04,934956,934956,4 to 34,2.131,0,0,0,0,0,0,0.000,0.927,0,0.000,414.511,180.549,0,0,0,0,0,315122,21.483,0.831,0,0,406.195,180.549,0,14,0,14,0,3343,9.475,0.791,0,0,408.542,97.363,0,15,0,15,0,719,5.683,0.464,0,0,404.546,180.549
pr05,1199566,1199566,4 to 36,2.154,0,0,0,0,0,0,0.000,0.927,0,0.000,389.001,189.796,0,0,0,0,0,318473,15.528,0.840,0,0,383.266,189.796,0,18,0,18,0,6569,7.054,0.803,0,0,382.554,85.237,0,14,1,13,0,1260,3.137,0.479,0,0,379.060,189.795
pr06,1682112,1682112,4 to 36,2.319,0,0,0,0,0,1,0.999,0.903,3,0.000,404.032,180.946,0,0,0,0,0,494874,15.421,0.884,0,0,397.125,180.946,0,13,0,13,0,6065,7.220,0.850,0,0,397.667,86.485,0,13,0,13,0,1032,2.169,0.479,0,0,394.081,180.946
pr07,149042,149042,4 to 28,2.075,0,0,0,0,0,0,0.000,0.762,0,0.000,357.194,188.866,0,0,0,0,0,39990,17.770,0.566,0,0,348.917,188.866,0,0,0,0,0,393,5.835,0.559,0,0,351.547,80.922,0,0,0,0,0,18,3.512,0.401,0,0,347.247,188.866
pr08,290587,290587,4 to 34,2.248,0,0,0,0,0,0,0.000,0.834,0,0.000,391.672,148.164,0,0,0,0,0,86704,10.062,0.801,0,0,385.396,148.164,0,3,0,3,0,2036,3.654,0.803,0,0,386.019,60.057,0,2,0,2,0,1774,2.534,0.498,0,0,381.852,148.164
pr09,400605,400605,4 to 36,2.218,0,0,0,0,0,2,1.342,0.804,0,0.000,396.569,138.729,0,0,0,0,0,146723,8.077,0.955,0,0,390.016,138.729,0,120,0,120,0,36057,2.323,0.936,0,0,389.730,55.898,0,222,0,222,0,19025,1.980,0.545,0,0,386.826,138.740
pr10,642027,642027,4 to 38,2.530,0,0,0,0,0,0,0.000,0.871,0,0.000,437.008,133.809,0,0,0,0,0,235329,8.103,1.088,0,0,429.214,133.809,0,83,0,83,0,15152,3.257,1.061,0,0,429.591,66.378,0,66,0,66,0,20754,2.242,0.601,0,0,427.197,133.807
pr11,219479,219479,4 to 26,2.012,0,0,0,0,0,0,0.000,0.567,0,0.000,356.472,272.798,0,0,0,0,0,44629,25.045,0.493,0,0,350.531,272.798,0,3,0,3,0,165,10.072,0.487,0,0,353.507,114.163,0,3,0,3,0,7,9.491,0.387,0,0,346.513,272.798
pr12,762104,762104,4 to 30,2.140,0,0,0,0,0,0,0.000,0.751,0,0.000,363.070,213.181,0,0,0,0,0,184521,27.090,0.772,0,0,362.406,213.181,0,9,0,9,0,547,6.062,0.750,0,0,356.877,75.641,0,7,0,7,0,105,3.642,0.489,0,0,353.251,213.180
pr13,1206019,1206019,4 to 32,2.234,0,0,0,0,0,0,0.000,0.687,0,0.000,369.229,159.462,0,0,0,0,0,427127,22.045,0.798,0,0,369.085,159.462,0,372,0,372,0,1935,8.359,0.769,0,0,367.712,58.110,0,341,3,338,0,313,4.360,0.514,0,0,359.873,159.448
pr14,1859557,1859557,4 to 34,2.545,0,0,0,0,0,0,0.000,0.715,0,0.000,364.596,153.785,0,0,0,0,0,382352,25.436,0.867,0,0,360.751,153.785,0,469,0,469,0,2242,6.640,0.830,0,0,358.839,50.754,0,410,15,395,0,1477,4.751,0.556,0,0,355.467,153.771
pr15,2965126,2965126,4 to 38,2.219,0,0,0,0,0,0,0.000,0.824,0,0.000,391.953,203.757,0,0,0,0,0,818035,27.374,0.960,0,0,390.874,203.757,0,196,0,196,0,10909,7.876,0.917,0,0,386.096,79.922,0,119,8,111,0,1644,7.601,0.532,0,0,382.086,203.756
pr16,4303839,4303839,4 to 42,2.541,0,0,0,0,0,24,2.460,0.855,9,0.000,389.023,165.978,0,1,0,1,0,1042878,24.480,0.977,0,0,387.105,165.978,0,132,0,132,0,6736,7.036,0.981,0,0,383.602,69.015,0,159,20,139,0,1161,5.260,0.570,0,0,379.427,165.978
pr17,328588,328588,4 to 28,1.902,0,0,0,0,0,0,0.000,0.541,0,0.000,336.010,197.822,0,0,0,0,0,67848,23.784,0.601,0,0,333.101,197.822,0,6,0,6,0,406,6.852,0.568,0,0,328.744,62.586,0,6,0,6,0,133,5.549,0.441,0,0,326.764,197.822
pr18,895460,895460,4 to 36,2.303,0,0,0,0,0,0,0.000,0.810,0,0.000,418.685,175.677,0,0,0,0,0,234296,18.485,0.987,0,0,418.378,175.677,0,38,0,38,0,1654,6.959,0.966,0,0,414.098,74.657,0,35,0,35,0,213,3.866,0.534,0,0,408.815,175.676
pr19,1580455,1580455,4 to 38,2.729,0,0,0,0,0,0,0.000,0.729,0,0.000,426.799,152.959,0,0,0,0,0,365403,12.652,1.184,0,0,424.179,152.959,0,19,0,19,0,4544,7.957,1.152,0,0,419.764,53.964,0,12,0,12,0,499,5.140,0.627,0,0,417.182,152.959
pr20,1627462,1627462,4 to 40,2.470,0,0,0,0,0,0,0.000,0.681,0,0.000,411.816,142.123,0,0,0,0,0,368965,14.941,1.165,0,0,409.276,142.123,0,61,0,61,0,14655,4.499,1.125,0,0,404.891,36.854,0,92,0,92,0,4668,3.275,0.621,0,0,402.582,142.120
\end{filecontents*}
\DTLloaddb{darp}{excel/darp/Table_logs.csv}
	\begin{table}[tb]
		\caption{Scheduling algorithms on DARP instances: number of incorrect infeasible declarations, number of deviating solutions, and average relative deviation.}
		\centering
		\label{table_darp1}
	    \begingroup
	    \setlength{\arrayrulewidth}{1.2pt}
		\resizebox{13.5cm}{!}{
			\begin{tabular}{lrr:rrr:rrr:rrr:rrr}\toprule
				\multicolumn{3}{c}{}
				& \multicolumn{3}{c}{Algorithm~(\ref{alg:algorithm1})}
				& \multicolumn{3}{c}{\cite{cordeaulaporte2003tabu}}
				& \multicolumn{3}{c}
                {\cite{parragh2009heuristic}}
                & \multicolumn{3}{c}
                {\cite{molenbruch2017multidir}}\\

				\cmidrule(lr){1-15}
				\bfseries Name & \bfseries \#Routes & \bfseries Size range & \bfseries \#Inf. & \bfseries \#Dev. & \bfseries Avg. dev. \% & \bfseries \#Inf. & \bfseries \#Dev. & \bfseries Avg. dev. \% & \bfseries \#Inf. & \bfseries \#Dev. & \bfseries Avg. dev. \% & \bfseries \#Inf. & \bfseries \#Dev. & \bfseries Avg. dev. \%
				\DTLforeach*{darp}{\Name=scenario,\numroutes=number of routes, \rangeroutes=range of route nodes,\cpuLP=average cpu LP [ms],\BONGincorr=BONG: number incorrect infeasible declarations,\BONGdev=BONG: number of deviating solutions, \BONGavgdev=BONG: average relative deviation, \BONcpu=BONG: average cpu algorithm [ms], \CORDincorr=CORD: number of incorrect infeasible declarations, \CORDdev=CORD: number of deviating solutions, \CORDavgdev=CORD: average relative deviation, \CORDcpu=CORD: average cpu algorithm [ms], \PARRincorr=PARR: number incorrect infeasible declarations, \PARRdev=PARR: number of deviating solutions, \PARRavgdev=PARR: average relative deviation, \PARRcpu=PARR: average cpu algorithm [ms], \MOLincorr=MOL: number incorrect infeasible declarations, \MOLdev=MOL: number of deviating solutions, \MOLavgdev=MOL: average relative deviation, \MOLcpu=MOL: average cpu algorithm [ms]}{
				\DTLiffirstrow{\\\hdashline}{\\}
					\Name & \numroutes & \rangeroutes & \BONGincorr & \BONGdev & \BONGavgdev & \CORDincorr & \CORDdev & \CORDavgdev & \PARRincorr & \PARRdev & \PARRavgdev & \MOLincorr & \MOLdev & \MOLavgdev}
				\\\bottomrule
			\end{tabular}}
	        \endgroup
	\end{table}

    Table~\ref{table_darp1} shows a comparison between the scheduling procedures on DARP instances. The first column indicates the instance name, following the convention employed in \cite{cordeaulaporte2003tabu}. These are instances that consider three to 13 vehicles and 24 to 144 requests. For further information on the generation of these DARP instances, the reader is referred to \cite{cordeaulaporte2003tabu}. The second column indicates the total number of evaluated feasible routes, and the third column the range of their size, i.e., the number of nodes in the routes. The fourth to sixt column show the total number of incorrect infeasible declarations, the total number of deviating solutions, and their absolute average percentage deviation from the optimal solution, by employing Algorithm~(\ref{alg:algorithm1}). The following columns show the same results for the 8-step scheduling algorithms by \cite{cordeaulaporte2003tabu}, \cite{parragh2009heuristic}, and \cite{molenbruch2017multidir}. As it can be noted, Algorithm~(\ref{alg:algorithm1}) does not return any incorrect infeasibility declaration. However, in some very rare cases (i.e., in only 27 occurrences over more than 21.5 million cases), it does return suboptimal solutions, namely for instances pr06,  pr09, and, mostly, pr16. For these instances, the returned suboptimal solutions deviate by less than 2.5\%, on absolute average terms. This consists of a very small proportion of suboptimal solutions, compared to the results obtained by employing the scheduling algorithms by  \cite{cordeaulaporte2003tabu}, \cite{parragh2009heuristic}, and \cite{molenbruch2017multidir}. Consistently with reported results from the literature, the approach by \cite{cordeaulaporte2003tabu} has the tendency to minimize the total number of incorrect infeasibility declarations at the cost of a higher number of suboptimal solutions and their average deviations. Namely, the approach returns about 27\% suboptimal solutions, with deviations between about 8\% and 27\%. The approach by \cite{parragh2009heuristic}, instead, has the tendency to minimize the total number of suboptimal solutions and their average deviation, at the cost of a higher number of incorrect infeasibility declarations. Namely, the approach returns up to 0.007\% incorrect infeasibility declarations but reduces the number of suboptimal solutions to about 0.5\%, with deviations between about 2.5\% and 10\%. The procedure by \cite{molenbruch2017multidir} produces similar results with respect to the procedure by \cite{parragh2009heuristic}, most specifically in terms of the number of incorrect infeasibility declarations. However, the the number of suboptimal solutions to 0.23\%, with deviations between about 2\% and 9.5\%. It is worth noting that, in our numerical experiments, all incorrect infeasibility declarations returned by the scheduling algorithms by \cite{parragh2009heuristic} and  \cite{molenbruch2017multidir} are due to violations of maximum ride time constraints. In addition to producing solutions of higher quality with respect to state-of-the-art scheduling algorithms, Algorithm~(\ref{alg:algorithm1}) is computationally efficient, as shown in Table~\ref{table_darp3}. It is about 60\% faster than a linear program and comparably efficient with respect to the scheduling algorithms by \cite{cordeaulaporte2003tabu}, \cite{parragh2009heuristic}, and \cite{molenbruch2017multidir}, with a slight computational disadvantage for small-scale instances, e.g. pr01, and advantage for large-scale instances, e.g., instances pr18, pr19, and pr20. 

 	\begin{table}[tb]
		\caption{Scheduling algorihtms on DARP instances: CPU times [ms]}
		\centering
		\label{table_darp3}
	    \begingroup
	    \setlength{\arrayrulewidth}{1.2pt}
		\resizebox{10cm}{!} {
			\begin{tabular}{l:r:r:r:r:r}\toprule
				\multicolumn{1}{c}{}
				& \multicolumn{1}{c}{Linear Program}
				& \multicolumn{1}{r}{Algorithm~(\ref{alg:algorithm1})}
				& \multicolumn{1}{c}{\cite{cordeaulaporte2003tabu}}
				& \multicolumn{1}{c}{\cite{parragh2009heuristic}}
                & \multicolumn{1}{c}
                {\cite{molenbruch2017multidir}}\\

				\cmidrule(lr){1-6}
				\bfseries Name & \bfseries Avg. CPU [ms] & \bfseries Avg. CPU [ms] & \bfseries Avg. CPU [ms] & \bfseries Avg. CPU [ms] & \bfseries Avg. CPU [ms]
				\DTLforeach*{darp}{\Name=scenario,\numroutes=number of routes, \rangeroutes=range of route nodes,\cpuLP=average cpu LP [ms],\BONGincorr=BONG: number incorrect infeasible declarations,\BONGdev=BONG: number of deviating solutions, \BONGavgdev=BONG: average relative deviation, \BONcpu=BONG: average cpu algorithm [ms], \CORDincorr=CORD: number of incorrect infeasible declarations, \CORDdev=CORD: number of deviating solutions, \CORDavgdev=CORD: average relative deviation, \CORDcpu=CORD: average cpu algorithm [ms], \PARRincorr=PARR: number incorrect infeasible declarations, \PARRdev=PARR: number of deviating solutions, \PARRavgdev=PARR: average relative deviation, \PARRcpu=PARR: average cpu algorithm [ms],\MOLincorr=MOL: number incorrect infeasible declarations, \MOLdev=MOL: number of deviating solutions, \MOLavgdev=MOL: average relative deviation, \MOLcpu=MOL: average cpu algorithm [ms]}{
				\DTLiffirstrow{\\\hdashline}{\\}
					\Name & \cpuLP & \BONcpu & \CORDcpu & \PARRcpu & \MOLcpu}
				\\\bottomrule
			\end{tabular}}
	        \endgroup
	\end{table} 

    \begin{filecontents*}{excel/uber/Table_logs.csv}
scenario,number of routes,number feasible routes,range of route nodes,average cpu LP [ms],BONG: number incorrect feasible declarations,BONG: number incorrect infeasible declarations,BONG: number infeasible declarations because of TW violations,BONG: number infeasible declarations because of RT violations,BONG: number infeasible declarations because of BATT violations,BONG: number of deviating solutions,BONG: average relative deviation,BONG: average cpu algorithm [ms],BONG: total number of recourse,BONG: avg recourse time,BONG: avg route completion time,BONG: total route wait time,CORD: number incorrect feasible declarations,CORD: number of incorrect infeasible declarations,CORD: number infeasible declarations because of TW violations,CORD: number infeasible declarations because of RT violations,CORD: number infeasible declarations because of BATT violations,CORD: number of deviating solutions,CORD: average relative deviation,CORD: average cpu algorithm [ms],CORD: total number of recourse,CORD: avg recourse time,CORD: avg route completion time,CORD: total route wait time,PARR: number incorrect feasible declarations,PARR: number incorrect infeasible declarations,PARR: number infeasible declarations because of TW violations,PARR: number infeasible declarations because of RT violations,PARR: number infeasible declarations because of BATT violations,PARR: number of deviating solutions,PARR: average relative deviation,PARR: average cpu algorithm [ms],PARR: total number of recourse,PARR: avg recourse time,PARR: avg route completion time,PARR: total route wait time,MOL: number incorrect feasible declarations,MOL: number incorrect infeasible declarations,MOL: number infeasible declarations because of TW violations,MOL: number infeasible declarations because of RT violations,MOL: number infeasible declarations because of BATT violations,MOL: number of deviating solutions,MOL: average relative deviation,MOL: average cpu algorithm [ms],MOL: total number of recourse,MOL: avg recourse time,MOL: avg route completion time,MOL: total route wait time
u2-16-0.1,37670,37670,4 to 23,1.945,0,0,0,0,0,0,0,0.676,0,0,110.600,64.026,0,0,0,0,0,14343,31.036,0.463,0,0,114.184,64.026,0,217,181,0,139,0,0.000,0.464,0,0,112.081,58.381,0,0,0,0,0,0,0.000,0.417,0,0,110.101,64.026
u2-16-0.4,20348,20348,5 to 23,1.878,0,0,0,0,0,0,0,0.663,0,0,112.365,63.362,0,0,0,0,0,3121,33.209,0.471,0,0,116.756,63.362,0,2460,154,0,2357,0,0.000,0.452,0,0,111.611,52.100,0,0,0,0,0,0,0.000,0.499,0,0,111.872,63.362
u2-16-0.7,15319,15319,6 to 24,1.966,0,0,0,0,0,0,0,0.503,0,0,107.714,40.019,0,0,0,0,0,1098,33.280,0.392,0,0,111.778,40.019,0,3736,63,0,3727,0,0.000,0.359,0,0,107.047,30.050,0,0,0,0,0,0,0.000,0.596,0,0,107.292,40.019
u2-20-0.1,51431,51431,5 to 30,2.062,0,0,0,0,0,0,0,0.910,0,0,137.716,76.866,0,166,0,166,0,12933,67.765,0.661,0,0,141.381,76.924,0,342,22,0,335,0,0.000,0.634,0,0,137.444,60.475,0,0,0,0,0,0,0.000,0.464,0,0,137.216,76.866
u2-20-0.4,39031,39031,5 to 31,2.655,0,0,0,0,0,0,0,1.025,0,0,145.990,95.039,0,13,0,13,0,9605,70.787,0.726,0,0,150.620,95.040,0,9144,24,0,9135,0,0.000,0.738,0,0,144.734,74.889,0,0,0,0,0,0,0.000,0.539,0,0,145.490,95.039
u2-20-0.7,31343,31343,5 to 27,2.099,0,0,0,0,0,0,0,0.774,0,0,147.760,79.442,0,7,0,7,0,5591,68.619,0.560,0,0,151.722,79.440,0,1296,11,0,1296,0,0.000,0.554,0,0,147.291,62.407,0,0,0,0,0,0,0.000,0.520,0,0,147.261,79.442
u2-24-0.1,26513,26513,6 to 33,2.872,0,0,0,0,0,0,0,1.034,0,0,172.664,79.168,0,0,0,0,0,10321,28.348,0.658,0,0,174.470,79.168,0,3609,1819,0,2596,0,0.000,0.649,0,0,172.204,69.819,0,0,0,0,0,0,0.000,0.629,0,0,172.168,79.168
u2-24-0.4,28913,28913,6 to 34,3.013,0,0,0,0,0,0,0,0.939,0,0,177.582,69.812,0,0,0,0,0,11061,26.981,0.657,0,0,179.183,69.813,0,3572,1065,0,3139,0,0.000,0.649,0,0,176.470,58.531,0,0,0,0,0,0,0.000,0.697,0,0,177.087,69.813
u2-24-0.7,1751,1751,5 to 27,2.556,0,0,0,0,0,0,0,1.108,0,0,112.182,100.041,0,0,0,0,0,86,24.093,0.976,0,0,119.719,100.041,0,315,42,0,312,0,0.000,0.941,0,0,111.506,26.745,0,0,0,0,0,0,0.000,1.771,0,0,111.797,100.041
u3-18-0.1,59678,59678,4 to 24,2.008,0,0,0,0,0,0,0,0.615,0,0,121.552,90.541,0,0,0,0,0,13580,36.428,0.377,0,0,127.224,90.541,0,1135,1135,0,152,0,0.000,0.365,0,0,124.757,73.592,0,0,0,0,0,0,0.000,0.355,0,0,121.058,90.541
u3-18-0.4,45877,45877,4 to 21,2.079,0,0,0,0,0,0,0,0.602,0,0,113.548,90.228,0,0,0,0,0,10629,38.876,0.345,0,0,116.735,90.228,0,384,357,0,364,0,0.000,0.348,0,0,118.438,70.468,0,0,0,0,0,0,0.000,0.353,0,0,113.056,90.228
u3-18-0.7,43889,43889,5 to 21,1.786,0,0,0,0,0,0,0,0.576,0,0,121.446,90.415,0,0,0,0,0,8149,46.895,0.373,0,0,126.215,90.415,0,7378,10,0,7377,0,0.000,0.356,0,0,122.719,63.570,0,0,0,0,0,0,0.000,0.383,0,0,120.961,90.415
u3-24-0.1,62295,62295,4 to 24,2.088,0,0,0,0,0,0,0,0.722,0,0,152.844,118.147,0,0,0,0,0,25272,59.241,0.459,0,0,156.141,118.147,0,38,36,0,14,7,0.565,0.468,0,0,155.992,106.356,0,0,0,0,0,0,0.000,0.389,0,0,152.346,118.147
u3-24-0.4,62863,62863,5 to 25,2.202,0,0,0,0,0,0,0,0.741,0,0,165.153,116.514,0,0,0,0,0,25325,58.769,0.495,0,0,168.498,116.514,0,3721,297,0,3692,2,0.565,0.496,0,0,166.360,104.547,0,0,0,0,0,0,0.000,0.434,0,0,164.657,116.514
u3-24-0.7,39789,39789,5 to 25,2.756,0,0,0,0,0,0,0,0.758,0,0,162.690,106.052,0,0,0,0,0,19159,50.038,0.529,0,0,166.978,106.052,0,4588,23,0,4584,1,0.565,0.515,0,0,163.559,91.079,0,0,0,0,0,0,0.000,0.488,0,0,162.196,106.052
u3-30-0.1,82775,82775,4 to 35,2.208,0,0,0,0,0,0,0,1.233,0,0,170.841,160.902,0,0,0,0,0,21557,106.719,0.758,0,0,177.728,160.902,0,75,55,0,52,0,0.000,0.740,0,0,175.665,123.086,0,0,0,0,0,0,0.000,0.453,0,0,170.343,160.902
u3-30-0.4,56560,56560,4 to 33,2.954,0,0,0,0,0,0,0,1.235,0,0,178.126,143.956,0,0,0,0,0,15463,63.967,0.761,0,0,184.130,143.956,0,8901,26,0,8885,0,0.000,0.750,0,0,185.649,113.062,0,0,0,0,0,0,0.000,0.540,0,0,177.628,143.956
u3-30-0.7,53357,53357,5 to 38,2.486,0,0,0,0,0,0,0,1.042,0,0,202.082,148.365,0,0,0,0,0,17616,65.755,0.623,0,0,206.915,148.366,0,20916,136,0,20863,0,0.000,0.572,0,0,207.924,123.280,0,0,0,0,0,0,0.000,0.551,0,0,201.590,148.365
u3-36-0.1,79254,79254,5 to 35,2.562,0,0,0,0,0,0,0,1.418,0,0,237.499,190.384,0,0,0,0,0,35159,68.424,0.912,0,0,241.938,190.385,0,6355,233,0,6779,0,0.000,0.885,0,0,235.158,156.836,0,0,0,0,0,0,0.000,0.526,0,0,237.000,190.384
u3-36-0.4,40374,40374,5 to 35,2.091,0,0,0,0,0,0,0,1.191,0,0,233.359,177.823,0,0,0,0,0,19270,72.027,0.755,0,0,239.432,177.823,0,7642,91,0,7610,0,0.000,0.704,0,0,230.246,139.520,0,0,0,0,0,0,0.000,0.505,0,0,232.862,177.823
u3-36-0.7,53003,53003,5 to 40,2.388,0,0,0,0,0,0,0,1.329,0,0,258.331,159.487,0,0,0,0,0,22234,59.382,0.829,0,0,262.127,159.487,0,7594,31,0,7594,0,0.000,0.836,0,0,261.558,141.701,0,0,0,0,0,0,0.000,0.616,0,0,257.832,159.487
u4-16-0.1,67699,67699,4 to 14,1.394,0,0,0,0,0,0,0,0.279,0,0,41.912,32.837,0,0,0,0,0,11342,96.120,0.194,0,0,46.142,32.837,0,0,0,0,0,0,0.000,0.193,0,0,43.367,15.407,0,0,0,0,0,0,0.000,0.247,0,0,41.415,32.837
u4-16-0.4,62234,62234,4 to 14,1.404,0,0,0,0,0,0,0,0.263,0,0,40.664,32.479,0,0,0,0,0,10092,101.253,0.187,0,0,44.575,32.479,0,0,0,0,0,0,0.000,0.185,0,0,42.248,15.167,0,0,0,0,0,0,0.000,0.253,0,0,40.173,32.479
u4-16-0.7,40835,40835,4 to 17,1.478,0,0,0,0,0,0,0,0.276,0,0,58.855,34.337,0,0,0,0,0,6270,62.277,0.249,0,0,63.368,34.337,0,5026,14,0,5026,0,0.000,0.247,0,0,56.915,10.162,0,0,0,0,0,0,0.000,0.337,0,0,58.456,34.337
u4-24-0.1,59968,59968,4 to 23,1.824,0,0,0,0,0,0,0,0.627,0,0,346.260,305.417,0,3,0,3,0,13711,42.533,0.329,0,0,350.116,305.425,0,334,304,0,360,0,0.000,0.328,0,0,353.256,298.366,0,0,0,0,0,0,0.000,0.333,0,0,345.668,305.417
u4-24-0.4,54756,54756,5 to 23,1.702,0,0,0,0,0,0,0,0.663,0,0,326.466,275.018,0,16,0,16,0,7858,30.573,0.379,0,0,329.675,275.063,0,8411,282,0,8385,0,0.000,0.366,0,0,315.018,243.644,0,0,0,0,0,0,0.000,0.373,0,0,325.881,275.018
u4-24-0.7,33597,33597,4 to 23,1.731,0,0,0,0,0,0,0,0.674,0,0,299.234,224.190,0,7,0,7,0,5010,34.655,0.417,0,0,302.122,224.216,0,2722,19,0,2721,0,0.000,0.414,0,0,295.088,203.463,0,0,0,0,0,0,0.000,0.415,0,0,298.651,224.190
u4-32-0.1,106225,106225,4 to 27,1.958,0,0,0,0,0,0,0,0.806,0,0,143.342,101.983,0,338,0,338,0,53488,48.561,0.522,0,0,149.078,101.954,0,754,662,0,483,0,0.000,0.513,0,0,148.737,87.140,0,0,0,0,0,0,0.000,0.415,0,0,142.846,101.983
u4-32-0.4,46931,46931,4 to 25,1.875,0,0,0,0,0,0,0,0.684,0,0,155.311,104.113,0,29,0,29,0,16534,54.786,0.500,0,0,158.964,104.116,0,8713,111,0,8686,0,0.000,0.485,0,0,157.045,91.100,0,0,0,0,0,0,0.000,0.428,0,0,154.820,104.113
u4-32-0.7,51616,51616,5 to 25,2.076,0,0,0,0,0,0,0,0.636,0,0,151.695,96.019,0,129,0,129,0,16893,44.910,0.472,0,0,156.885,96.010,0,14747,38,0,14747,0,0.000,0.497,0,0,150.363,66.876,0,0,0,0,0,0,0.000,0.421,0,0,151.208,96.019
u4-40-0.1,46790,46790,4 to 29,2.216,0,0,0,0,0,0,0,1.064,0,0,366.539,378.017,0,0,0,0,0,20143,67.739,0.643,0,0,376.710,378.017,0,4635,184,0,4772,0,0.000,0.624,0,0,367.346,305.556,0,0,0,0,0,0,0.000,0.486,0,0,366.042,378.016
u4-40-0.4,44843,44843,5 to 31,2.031,0,0,0,0,0,0,0,0.918,0,0,392.990,353.199,0,0,0,0,0,22624,86.456,0.541,0,0,399.849,353.199,0,13034,52,0,13032,0,0.000,0.535,0,0,394.897,304.694,0,0,0,0,0,0,0.000,0.488,0,0,392.494,353.199
u4-40-0.7,31625,31625,5 to 25,1.490,0,0,0,0,0,0,0,0.369,0,0,277.830,305.410,0,0,0,0,0,1434,59.386,0.273,0,0,283.990,305.410,0,3516,364,0,3509,0,0.000,0.261,0,0,273.167,219.130,0,0,0,0,0,0,0.000,0.383,0,0,277.476,305.410
u4-48-0.1,90666,90666,5 to 39,2.445,0,0,0,0,0,0,0,1.339,0,0,240.387,177.752,0,175,0,175,0,57814,37.664,0.990,0,0,243.677,177.722,0,9379,221,1,12731,0,0.000,0.964,0,0,240.815,154.353,0,1,0,1,0,0,0.000,0.543,0,0,239.889,177.752
u4-48-0.4,58420,58420,5 to 31,1.647,0,0,0,0,0,0,0,0.407,0,0,159.147,145.273,0,107,0,107,0,5979,48.564,0.345,0,0,163.673,145.155,0,2731,1448,1,2497,0,0.000,0.339,0,0,162.162,100.685,0,1,0,1,0,0,0.000,0.400,0,0,158.828,145.273
u4-48-0.7,43857,43857,5 to 29,1.605,0,0,0,0,0,0,0,0.411,0,0,156.159,136.618,0,115,0,115,0,4301,44.029,0.340,0,0,160.763,136.413,0,3524,139,1,3523,0,0.000,0.345,0,0,157.227,90.308,0,1,0,1,0,0,0.000,0.405,0,0,155.801,136.617
u5-40-0.1,130743,130743,4 to 25,1.945,0,0,0,0,0,0,0,0.849,0,0,149.655,126.202,0,0,0,0,0,59482,77.189,0.515,0,0,156.118,126.202,0,237,82,0,213,0,0.000,0.494,0,0,152.849,91.998,0,0,0,0,0,0,0.000,0.390,0,0,149.062,126.202
u5-40-0.4,116848,116848,4 to 27,2.145,0,0,0,0,0,0,0,0.795,0,0,169.808,120.363,0,0,0,0,0,44558,79.516,0.543,0,0,174.531,120.363,0,6327,227,0,6289,0,0.000,0.532,0,0,171.049,93.711,0,0,0,0,0,0,0.000,0.436,0,0,169.214,120.363
u5-40-0.7,50576,50576,5 to 23,1.325,0,0,0,0,0,0,0,0.204,0,0,89.798,93.350,0,0,0,0,0,2407,63.427,0.199,0,0,92.643,93.350,0,3962,811,0,3962,1,7.262,0.195,0,0,89.232,46.720,0,0,0,0,0,0,0.000,0.299,0,0,89.544,93.350
u5-50-0.1,112835,112835,4 to 34,2.185,0,0,0,0,0,0,0,1.057,0,0,193.501,179.413,0,1,0,1,0,47249,60.896,0.645,0,0,201.735,179.413,0,4029,93,0,4120,2,3.583,0.629,0,0,196.157,145.192,0,0,0,0,0,0,0.000,0.429,0,0,193.003,179.413
u5-50-0.4,111003,111003,4 to 31,2.208,0,0,0,0,0,0,0,1.043,0,0,222.417,168.935,0,3,0,3,0,61626,60.273,0.692,0,0,228.931,168.934,0,27094,60,0,27077,6,5.335,0.678,0,0,222.589,133.000,0,0,0,0,0,1,3.041,0.461,0,0,221.922,168.935
u5-50-0.7,76143,76143,5 to 33,2.306,0,0,0,0,0,0,0,1.065,0,0,218.598,164.566,0,0,0,0,0,41223,52.761,0.702,0,0,226.739,164.566,0,10444,3,0,10444,123,6.021,0.721,0,0,215.465,118.798,0,0,0,0,0,4,4.070,0.507,0,0,218.103,164.566
\end{filecontents*}
\DTLloaddb{uber}{excel/uber/Table_logs.csv}
	\begin{table}[tb]
		\caption{Scheduling algorithms on E-ADARP instances: number of incorrect infeasible declarations, number of deviating solutions, and average relative deviation.}
		\centering
		\label{table_uber1}
	    \begingroup
	    \setlength{\arrayrulewidth}{1.2pt}
		\resizebox{12cm}{!} {
			\begin{tabular}{lrr:rr:rrr:rrr:rrr}\toprule
				\multicolumn{3}{c}{}
				& \multicolumn{2}{c}{Algorithm~(\ref{alg:algorithm1})}
				& \multicolumn{3}{c}{\cite{cordeaulaporte2003tabu}}
				& \multicolumn{3}{c}{\cite{parragh2009heuristic}}
                & \multicolumn{3}{c}
                {\cite{molenbruch2017multidir}}\\

				\cmidrule(lr){1-14}
				\bfseries Name & \bfseries \#Routes & \bfseries Size range & \bfseries \#Inf. & \bfseries \#Dev. & \bfseries \#Inf. & \bfseries \#Dev. & \bfseries Avg. dev. \% & \bfseries \#Inf. & \bfseries \#Dev. & \bfseries Avg. dev. \% & \bfseries \#Inf. & \bfseries \#Dev. & \bfseries Avg. dev. \% 
				\DTLforeach*{uber}{\Name=scenario,\numroutes=number of routes, \rangeroutes=range of route nodes,\cpuLP=average cpu LP [ms],\BONGincorr=BONG: number incorrect infeasible declarations,\BONGdev=BONG: number of deviating solutions, \BONGavgdev=BONG: average relative deviation, \BONcpu=BONG: average cpu algorithm [ms], \CORDincorr=CORD: number of incorrect infeasible declarations, \CORDdev=CORD: number of deviating solutions, \CORDavgdev=CORD: average relative deviation, \CORDcpu=CORD: average cpu algorithm [ms], \PARRincorr=PARR: number incorrect infeasible declarations, \PARRdev=PARR: number of deviating solutions, \PARRavgdev=PARR: average relative deviation, \PARRcpu=PARR: average cpu algorithm [ms],\MOLincorr=MOL: number incorrect infeasible declarations, \MOLdev=MOL: number of deviating solutions, \MOLavgdev=MOL: average relative deviation, \MOLcpu=MOL: average cpu algorithm [ms]}{
				\DTLiffirstrow{\\\hdashline}{\\}
					\Name & \numroutes & \rangeroutes & \BONGincorr & \BONGdev & \CORDincorr & \CORDdev & \CORDavgdev & \PARRincorr & \PARRdev & \PARRavgdev & \MOLincorr & \MOLdev & \MOLavgdev}
				\\\bottomrule
			\end{tabular}}
	        \endgroup
	\end{table} 

    Table~\ref{table_uber1} shows a comparison between the investigated scheduling procedures on e-ADARP instances. The first column indicates the instance name, following the convention employed in \cite{bongiovanni2019}, i.e., $<$u$>$$<$number of vehicles$>$-$<$number of customers$>$-$<$minimum battery inventory ratio at the destination depot$>$,``u" is used to refer to the instances adapted from real data from Uber Technologies Inc. Furthermore, the reader is reminded that the minimum battery inventory ratio is used to compute the minimal battery levels that all vehicles need to have at the destination depot, i.e. a minimal battery ratio of 0.7 means vehicles need to have at least 70\% of their nominal battery capacity at the destination depot. The following columns adopt the same convention used in Table~\ref{table_darp1}. As it can be noted, in this case Algorithm~(\ref{alg:algorithm1}) always returns optimal scheduling solutions, which include charging decisions. Differently from the DARP results shown in Table~\ref{table_darp1}, the approach of \cite{cordeaulaporte2003tabu} return some incorrect infeasibility declarations, although this remains limited to about 0.05\%. The number of suboptimal solutions, however, rise up to about 34\% and so their absolute average deviations from the optimal solution, which range between about 24\% and more than 100\%. The approach of \cite{parragh2009heuristic} also shows an increase in the total number of incorrect infeasibility declarations with respect to the DARP results shown in Table~\ref{table_darp1}. Namely, the total number of infeasibility declarations rise to at most about 10\%, however, the number of suboptimal solutions decrease to about 0.006\%, with deviations which are contained between about 0.5\% and 7\%. Finally, the approach of \cite{molenbruch2017multidir} shows a net decrease both in terms of the total number of incorrect infeasibility declarations and suboptimal solutions, with a total of eight occurrences over more than 2 million cases, with respect to the DARP results shown in Table~\ref{table_darp1}.

	  \begin{table}[tb]
		\caption{Scheduling algorithms on E-ADARP instances: infeasibility reasons}
		\centering
		\label{table_uber2}
	    \begingroup
	    \setlength{\arrayrulewidth}{1.2pt}
		\resizebox{9cm}{!} {
			\begin{tabular}{l:rrr:rrr:rrr}\toprule
				\multicolumn{1}{c}{}
				& \multicolumn{3}{c}{\cite{cordeaulaporte2003tabu}}
				& \multicolumn{3}{c}{\cite{parragh2009heuristic}}
                & \multicolumn{3}{c}
                {\cite{molenbruch2017multidir}}\\

				\cmidrule(lr){1-10}
				\bfseries Name & \bfseries \#Inf. TW & \bfseries \#Inf. RT & \bfseries \#Inf. BATT & \bfseries \#Inf. TW & \bfseries \#Inf. RT & \bfseries \#Inf. BATT & \bfseries \#Inf. TW & \bfseries \#Inf. RT & \bfseries \#Inf. BATT
				\DTLforeach*{uber}{\Name=scenario,\numroutes=number of routes, \rangeroutes=range of route nodes,\cpuLP=average cpu LP [ms],\BONGincorr=BONG: number incorrect infeasible declarations,\BONGdev=BONG: number of deviating solutions, \BONGavgdev=BONG: average relative deviation, \BONcpu=BONG: average cpu algorithm [ms], \CORDincorr=CORD: number of incorrect infeasible declarations, \CORDtw=CORD: number infeasible declarations because of TW violations, \CORDrt=CORD: number infeasible declarations because of RT violations, \CORDbatt=CORD: number infeasible declarations because of BATT violations, \CORDdev=CORD: number of deviating solutions, \CORDavgdev=CORD: average relative deviation, \CORDcpu=CORD: average cpu algorithm [ms], \PARRincorr=PARR: number incorrect infeasible declarations, \PARRtw=PARR: number infeasible declarations because of TW violations, \PARRrt=PARR: number infeasible declarations because of RT violations, \PARRbatt=PARR: number infeasible declarations because of BATT violations, \PARRdev=PARR: number of deviating solutions, \PARRavgdev=PARR: average relative deviation, \PARRcpu=PARR: average cpu algorithm [ms],\MOLincorr=MOL: number incorrect infeasible declarations, \MOLtw=MOL: number infeasible declarations because of TW violations, \MOLrt=MOL: number infeasible declarations because of RT violations, \MOLbatt=MOL: number infeasible declarations because of BATT violations, \MOLdev=MOL: number of deviating solutions, \MOLavgdev=MOL: average relative deviation, \MOLcpu=MOL: average cpu algorithm [ms]}{
				\DTLiffirstrow{\\\hdashline}{\\}
					\Name & \CORDtw & \CORDrt & \CORDbatt & \PARRtw & \PARRrt & \PARRbatt & \MOLtw & \MOLrt & \MOLbatt}
				\\\bottomrule
			\end{tabular}}
	        \endgroup
	\end{table} 
	
	\begin{table}[tbh]
		\caption{Scheduling algorihtms on E-ADARP instances: CPU times [ms]}
		\centering
		\label{table_uber3}
	    \begingroup
	    \setlength{\arrayrulewidth}{1.2pt}
		\resizebox{11cm}{!} {
			\begin{tabular}{l:r:r:r:r:r}\toprule
				\multicolumn{1}{c}{}
				& \multicolumn{1}{c}{Linear Program}
				& \multicolumn{1}{c}{Algorithm~(\ref{alg:algorithm1})}
				& \multicolumn{1}{c}{\cite{cordeaulaporte2003tabu}}
				& \multicolumn{1}{c}{\cite{parragh2009heuristic}}
                & \multicolumn{1}{c}
                {\cite{molenbruch2017multidir}}\\

				\cmidrule(lr){1-6}
				\bfseries Name & \bfseries Avg. CPU [ms] & \bfseries Avg. CPU [ms] & \bfseries Avg. CPU [ms] & \bfseries Avg. CPU [ms] & \bfseries Avg. CPU [ms]
				\DTLforeach*{uber}{\Name=scenario,\numroutes=number of routes, \rangeroutes=range of route nodes,\cpuLP=average cpu LP [ms],\BONGincorr=BONG: number incorrect infeasible declarations,\BONGdev=BONG: number of deviating solutions, \BONGavgdev=BONG: average relative deviation, \BONcpu=BONG: average cpu algorithm [ms], \CORDincorr=CORD: number of incorrect infeasible declarations, \CORDdev=CORD: number of deviating solutions, \CORDavgdev=CORD: average relative deviation, \CORDcpu=CORD: average cpu algorithm [ms], \PARRincorr=PARR: number incorrect infeasible declarations, \PARRdev=PARR: number of deviating solutions, \PARRavgdev=PARR: average relative deviation, \PARRcpu=PARR: average cpu algorithm [ms],\MOLincorr=MOL: number incorrect infeasible declarations, \MOLdev=MOL: number of deviating solutions, \MOLavgdev=MOL: average relative deviation, \MOLcpu=MOL: average cpu algorithm [ms]}{
				\DTLiffirstrow{\\\hdashline}{\\}
					\Name & \cpuLP & \BONcpu & \CORDcpu & \PARRcpu & \MOLcpu}
				\\\bottomrule
			\end{tabular}}
	        \endgroup
	\end{table}

    As shown in Table~\ref{table_uber2}, in the case of an e-ADARP, the procedure by \cite{cordeaulaporte2003tabu} mostly produce solutions that violate ride time constraints, whereas the procedure by \cite{parragh2009heuristic} mostly violate time window and battery constraints. This last procedure is in fact designed to minimize completion time, which, in turn, may increase the chances of violating battery management constraints. As for the DARP results shown in Table~\ref{table_darp1}, the procedure by \cite{molenbruch2017multidir} may violate ride time constraints. Finally, Table~\ref{table_uber3} shows the average CPU times, in milliseconds, obtained by running the linear program, Algorithm~(\ref{alg:algorithm1}), and the scheduling procedures by \cite{cordeaulaporte2003tabu}, \cite{parragh2009heuristic}, and \cite{molenbruch2017multidir} on the e-ADARP. As it can be noted, also in this case,  Algorithm~(\ref{alg:algorithm1}) is about 60\% faster with respect to a linear program, on average, and up to about 80-85\% faster, e.g. for instance u5-40-0.7 and instances u4-16. Savings in terms of computing time might be even more significant when considering larger problem instances or when several static problems are solved in real time, e.g. as in the dynamic e-ADARP presented in \cite{bongiovanni2022machine}.

\section{Summary}\label{conclusion}

This work proposed an excess-ride-time procedure for scheduling dial-a-ride instances. The procedure is shown to produce high-quality solutions, which are only very rarely suboptimal, in at least half the time of a linear program and in comparable time with respect to state-of-the-art scheduling algorithms. The proposed procedure can be potentially employed to efficiently solve scheduling problems from static and dynamic metaheuristic appraches. For the e-ADARP, we further propose a battery heuristic which can be used when the vehicle routes include one or more visits to charging facilities. The heuristic is needed to provide battery-feasible charging plans for the vehicle schedules. The battery heuristic builds on the assumption that charging as much as possible as early as possible is the best strategy for battery-feasibility.

Computational experiments are carried on a plethora of static DARP and e-ADARP instances from the literature, which are extracted from an adaptive large neighborhood search with 1,000 iterations. Experiments show that the route evaluation procedure is computationally more efficient than solving a linear program, while returning higher-quality solutions with respect to popular scheduling heuristics from the literature.  Furthermore, results show that the proposed scheduling procedure does not produce incorrect infeasibility declarations and produces suboptimal solutions in very rare cases (i.e., in about 0.0001 \% of times, which deviate up to 2.5\% from optimal solutions, on average).

\section*{Acknowledments}
\noindent
This research was enabled in part by support provided by Calcul Qu\'{e}bec (\url{www.calculquebec.ca/en/}) and Compute Canada (\url{www.computecanada.ca}).

\bibliographystyle{elsarticle-harv}


\end{document}